**An approximation scheme for optimal control of Volterra integral equations**

S. A. Belbas
Mathematics Department
University of Alabama
Tuscaloosa, AL. 35487-0350. USA.

e-mail: SBELBAS@GP.AS.UA.EDU

Abstract. We present and analyze a new method for solving optimal control problems for Volterra integral equations, based on approximating the controlled Volterra integral equations by a sequence of systems of controlled ordinary differential equations. The resulting approximating problems can then be solved by dynamic programming methods for ODE controlled systems. Other, straightforward versions of dynamic programming, are not applicable to Volterra integral equations. We also derive the connection between our version of dynamic programming and the Hamiltonian equations for Volterra controlled systems.

1. Introduction and statement of the problem.

We are interested in a controlled system governed by a Volterra integral equation

$$x(t) = x_0(t) + \int_0^t f(t, s, x(s), u(s))\, ds$$

--- (1.1)

and the associated optimal control problem of minimizing a cost functional of the form

$$J := \int_0^T F(t, x(t), u(t))\, dt + F_0(x(T))$$

--- (1.2)

In this paper, we develop and analyze a novel method for solving this optimal control problem. Our method consists of reducing the original control problem to an equivalent problem for an infinite system of ordinary differential equations, then approximating the infinite-dimensional problem by a truncated finite-dimensional problem, solving that finite dimensional approximation by the method of dynamic programming for controlled ordinary differential equations, and then passing to the limit as the dimension of the finite-dimensional truncation goes to infinity. The theory of optimal control of ordinary differential equations uses two main groups of methods: methods based on dynamic programming, and necessary conditions of the type of Pontryagin's maximum principle. Problems of optimal control for Volterra integral equations have been treated by the method of necessary conditions of a type akin to Pontryagin's maximum principle.

The method of dynamic programming has not found, up to now, applications in the area of optimal control of Volterra integral equations. This has been due to the lack of a suitable analytical framework for the applicability of dynamic programming techniques.

A naive attempt to utilize dynamic programming for controlled Volterra integral equations would utilize a parametrization of the state dynamics and the cost functional by time t, the part of the trajectory up to time t, and the history of the control function up to time t. With such a parametrization, the controlled Volterra equation becomes a dynamical system over an infinite-dimensional space.
Let the original controlled Volterra equation be

$$x(t) = x_0(t) + \int_0^t f(t, s, x(s), u(s))\, ds$$

--- (1.1)

We may define a new state $\tilde{x}(t)$ and a new control $\tilde{u}(t)$ as functions from [0, T] into $C([0,1] \mapsto \mathbb{R})$, as follows:

$$(\tilde{x}(t))(\tau) := x(t\tau), \ (\tilde{u}(t))(\tau) := u(t\tau); \ 0 \leq \tau \leq 1$$

--- (1.2)

Then, for given $\xi := \tilde{x}(t)$, $\beta := \tilde{u}(t)$, $\alpha(.) := u|_{(t,t+\delta t]}$, we can find $x(.)|_{(t,t+\delta t]}$ from the Volterra equation (1.1). In this way, the original Volterra equation becomes a dynamical system with state-space $C([0,1] \mapsto \mathbb{R})$. Once this dynamical system has been established, it is possible to apply the techniques of dynamic programming, under suitable assumptions. However, this type of dynamic programming offers little advantage, as the state-space for the dynamic programming equations has the same dimensionality as the full trajectory $\{x(t): 0 \leq t \leq T\}$ obtained from the original Volterra equation (1.1).

## 2. Derivation of the infinite-dimensional dynamic programming equations.

We start with the system dynamics (1.1). With successive differentiations, formally at this stage and to be rigorously justified later, we find

$$x(t) = x_0(t) + x_1(t), \quad x_1(t) := \int_0^t f(t,s,x_0(s)+x_1(s),u(s))\,ds;$$

$$\dot{x}_1(t) = f(t,s,x_0(s)+x_1(s),u(s))|_{s=t} + x_2(t), \quad x_2(t) := \int_0^t \frac{\partial}{\partial t} f(t,s,x_0(s)+x_1(s),u(s))\,ds;$$

...
...

$$\dot{x}_n(t) = \frac{\partial^{n-1}}{\partial t^{n-1}} f(t,s,x_0(s)+x_1(s),u(s))|_{s=t} + x_{n+1}(t),$$

$$x_{n+1}(t) := \int_0^t \frac{\partial^n}{\partial t^n} f(t,s,x_0(s)+x_1(s),u(s))\,ds$$

...

--- (2.1)

We set

$$g_n(t,x,u) := \frac{\partial^{n-1}}{\partial t^{n-1}} f(t,s,x,u)|_{s=t}$$

--- (2.2)

Then (2.1) becomes a Cauchy problem for an infinite-dimensional differential system

$$\dot{x}_n(t) = g_n(t,x_0(t)+x_1(t),u(t)) + x_{n+1}(t), \quad x_n(0)=0;\ n=1,2,...$$

--- (2.3)

If $x(t)$ is the solution of (1.1), and if we define $x_n(t)$ by

$$x_n(t) := \int_0^t \frac{\partial^{n-1}}{\partial t^{n-1}} f(t, s, x_0(s) + x_1(s), u(s)) \, ds, \text{ for } n \geq 2;$$

$$x_1(t) := \int_0^t f(t, s, x_0(s) + x_1(s), u(s)) \, ds$$

--- (2.4)

then the collection $\{x_n(t) : n = 1, 2, ...\}$ solves the Cauchy problem (2.3). It follows that, if the solution $\{x_n(t) : n = 1, 2, ...\}$ of (2.3) is unique, then, by defining $x(t) := x_0(t) + x_1(t)$, we can conclude that $x(t)$ is a solution of (1.1).

Now, the problem of minimizing the functional J, defined in (1.2), under the state dynamics (1.1), becomes equivalent to minimizing J subject to the infinite-dimensional ODE system (2.3). We interpret (2.3) as an initial value problem for an ODE with state space $\ell^\infty$, the space of all bounded real valued sequences equipped with the supremum norm. The standard method of dynamic programming then leads to the following equation:

$$\inf_{\alpha \in K} \left\{ \frac{\partial V(t, \tilde{x})}{\partial t} + \sum_{i=1}^\infty \frac{\partial V(t, \tilde{x})}{\partial \tilde{x}_i} [g_i(t, x_0(t) + \tilde{x}_1, \alpha) + \tilde{x}_{i+1}] + F(t, x_0(t) + \tilde{x}_1, \alpha) \right\} = 0;$$

$$V(T, \tilde{x}) = F_0(x_0(T) + \tilde{x}_1)$$

--- (2.5)

In this paper, we shall not pursue the question of solving Volterra integral equations via reduction to infinite-dimensional ODE systems; we mention, however, that (i) this reduction is, to the best of our knowledge, an original result presented for the first time in this paper, and (ii) the theory of infinite-dimensional ODE systems is part of classical Mathematics, it is contained in [C] and in several other texts.

An important particular case arises when f is a polynomial in t, say

$$f(t, s, x, u) = \sum_{j=0}^N t^j f_j(s, x, u)$$

--- (2.6)

In that case, $\dfrac{\partial^m}{\partial t^m} f(t,s,x,u) \equiv 0$ for $m \geq N+1$, thus $x_m(t) \equiv 0$ for $m \geq N+2$, and then (2.3) becomes an (N+1)-dimensional ODE system, to which the standard theory of weak solutions of (finite-dimensional) Hamilton-Jacobi equations can be applied.

These considerations lead us to examine a family of polynomial (in t) approximations to the kernel f(t, s, x, u) of (1.1), and the corresponding approximations to the state x(t) and to the cost functional J.

A polynomial approximation to f is an expression of the form

$$f_{[N]}(t,s,x,u) = \sum_{j=0}^{N} t^j f_j(s,x,u)$$

--- (2.7)

The approximation is meaningful if, for every natural number N, the (N+1)-th derivative of f with respect to t is continuous in t and satisfies

$$\left| \dfrac{\partial^{N+1}}{\partial t^{N+1}} f(t,s,x,u) \right| \leq M_N = \text{const. for all relevant values of t, s, x, u.}$$

Then we can write

$$f(t,s,x,u) = f_{[N]}(t,s,x,u) + R_N(t,s,x,u)$$

--- (2.8)

where the remainder $R_N(t,s,x,u)$ satisfies

$$|R_N(t,s,x,u)| \leq \dfrac{M_N t^{N+1}}{(N+1)!} \text{ for all relevant values of t, s, x, u}$$

--- (2.9)

Let x(t) be the solution of (1.1) for some (arbitrary but fixed) choice of the control function u, and let $x_{[N]}(t)$ be the solution of (1.1) for the same choice of the control function u.
We assume that f is Lipschitz in x, uniformly in the other variables that appear in f:

$$|f(t,s,x,u) - f(t,s,y,u)| \leq L_0 |x-y| \text{ for all relevant values of t, s, x, y, u}$$

for some constant $L_0$.

The functions $x(t)$ and $x_{[N]}(t)$ are defined as the solutions of

$$x(t) = x_0(t) + \int_0^t f(t,s,x(s),u(s))ds, \quad x_{[N]}(t) = x_0(t) + \int_0^t f_{[N]}(t,s,x_{[N]}(s),u(s))ds$$

--- (2.11)

We shall prove:

<u>Theorem 2.1.</u> Under the stated assumptions, the sequence $\{x_{[N]}(\cdot): N = 1,2,...\}$ converges to $x(\cdot)$ uniformly for $t \in [0,T]$, as $N \to \infty$.

<u>Proof:</u> We have

$$|x(t) - x_{[N]}(t)| \leq \int_0^t |f(t,s,x(s),u(s)) - f_{[N]}(t,s,x_{[N]}(s),u(s))|\,ds \leq$$

$$\leq \int_0^t |f(t,s,x(s),u(s)) - f(t,s,x_{[N]}(s),u(s))|\,ds +$$

$$+ \int_0^t |f(t,s,x_{[N]}(s),u(s)) - f_{[N]}(t,s,x_{[N]}(s),u(s))|\,ds \leq$$

$$\leq \int_0^t L_0 |x(s) - x_{[N]}(s)|\,ds + \int_0^t M_N \frac{s^{N+1}}{(N+1)!}\,ds =$$

$$= \int_0^t L_0 |x(s) - x_{[N]}(s)|\,ds + M_N \frac{t^{N+2}}{(N+2)!}$$

--- (2.12)

By Gronwall's inequality, we have

$$|x(t) - x_{[N]}(t)| \leq z_N(t)$$

--- (2.13)

where $z_N(t)$ solves

$$z_N(t) = M_N \frac{t^{N+2}}{(N+2)!} + L_0 \int_0^t z_N(s)\,ds$$

--- (2.14)

By elementary methods, we find that the solution of (2.14) is

$$z_N(t) = M_N \left[ \frac{\exp(L_0 t) - 1}{L_0^{N+2}} - \sum_{\ell=0}^{N} \frac{t^{N-\ell+1}}{L_0^{\ell+1}[(N-\ell+1)!]} \right]$$

--- (2.15)

An equivalent expression for $z_N(t)$ is found by expanding $\exp(L_0 t)$ as a power series, thus obtaining

$$z_N(t) = \frac{M_N}{L_0^{N+2}} \sum_{k=N+2}^{\infty} \frac{(L_0 t)^k}{k!}$$

--- (2.16)

If the sequence $\{M_N : N = 1, 2, 3, ...\}$ is bounded, then (2.16) implies that $\lim_{N \to \infty} z_N(t) = 0$ uniformly for $t \in [0, T]$. Indeed, it follows from (2.16) that

$$0 \leq z_N(t) \leq z_N(T) = \frac{M_N}{L_0^{N+2}} \sum_{k=N+2}^{\infty} \frac{(L_0 T)^k}{k!}, \quad \forall t \in [0, T]$$

--- (2.17).

If $L_0 \leq 1$, then $\dfrac{M_N}{L_0^{N+2}} \sum_{k=N+2}^{\infty} \dfrac{(L_0 T)^k}{k!} \leq M_N \sum_{k=N+2}^{\infty} \dfrac{T^k}{k!}$, and $\sum_{k=N+2}^{\infty} \dfrac{T^k}{k!} \to 0$ as $N \to \infty$, since $\sum_{k=N+2}^{\infty} \dfrac{T^k}{k!}$ is the tail part of the power series for $\exp T$. If $L_0 > 1$, then

$$\frac{M_N}{L_0^{N+2}} \sum_{k=N+2}^{\infty} \frac{(L_0 T)^k}{k!} \leq \frac{M_N}{L_0^{N+2}} \sum_{k=0}^{\infty} \frac{(L_0 T)^k}{k!} = \frac{M_N \exp(L_0 T)}{L_0^{N+2}} \to 0 \text{ as } N \to \infty.$$

Thus in each of the two cases ($L_0 > 1$ or $L_0 \leq 1$) we conclude that $z_N(t) \to 0$ as $N \to \infty$ uniformly in $t \in [0, T]$. ///

Remark 2.1. The quantity $z_N(T)$ provides an estimate for the rate of convergence in theorem 2.1. ///

Now we examine the behaviour of the functional J under the polynomial approximation in the state dynamics. We assume, in addition to the previous conditions, that F and $F_0$ are Lipschitz in x, uniformly in the other variables, with Lipschitz constants $L_F$, $L_{F_0}$, respectively. We denote by $J(u(.))$ the value of the cost functional (1.2) obtained by using the control function u(.), and we denote by $J_{[N]}(u(.))$ the functional

$$J_{[N]}(u(.)) := \int_0^T F(t, x_{[N]}(t), u(t)) dt + F_0(x_{[N]}(T))$$

--- (2.18)

where $x_{[N]}(.)$ is solution of the second equation in (2.11). Then we have:

Theorem 2.2. The sequence of functionals $\{J_N : N = 1, 2, ...\}$ converges to J, as $N \to \infty$, uniformly with respect to u, and the rate of convergence is estimated by the quantity $\omega_N(T)$ given in (2.19) below.

Proof: We estimate the difference $|J(u(.)) - J_{[N]}(u(.))|$ as follows:

$$|J(u(.)) - J_{[N]}(u(.))| \leq L_F \int_0^T z_N(t) dt + L_{F_0} z_N(T) =$$

$$= L_F M_N \left[ \frac{\exp(L_0 T) - 1 - L_0 T}{L_0^{N+3}} - \sum_{\ell=0}^{N} \frac{T^{N-\ell+2}}{L_0^{\ell+1}[(N-\ell+2)!]} \right] +$$

$$+ L_{F_0} M_N \left[ \frac{\exp(L_0 T) - 1}{L_0^{N+2}} - \sum_{\ell=0}^{N} \frac{T^{N-\ell+1}}{L_0^{\ell+1}[(N-\ell+1)!]} \right] \equiv \omega_N(T)$$

--- (2.19)

It follows (by the same type of argument as in the proof of theorem 2.1) that $\omega_N(T) \to 0$ as $N \to \infty$. ///

Let

$$J^*_{[N]} := \inf_{u(.)\in U_{ad}} J_{[N]}(u(.)), \quad J^* := \inf_{u(.)\in U_{ad}} J(u(.))$$

--- (2.20)

and let $u^*_{\varepsilon,N}(.) \in U_{ad}$ be an ε-suboptimal control with respect to $J_{[N]}$, i.e.

$$J_{[N]}(u^*_{\varepsilon,N}(.)) \leq J^*_{[N]} + \varepsilon$$

--- (2.21)

Next, we prove:

Theorem 2.3. The following inequality holds:

$$J(u^*_{\varepsilon,N}(.)) - \varepsilon - 2\omega_N(T) \leq J^*(u(.)) \leq J(u^*_{\varepsilon,N}(.))$$

and, as a consequence, the sequence $\{J(u^*_{\varepsilon,N}(.)): N = 1,2,...\}$ converges to $J^*(u(.))$ as $N \to \infty$.

Proof: One side of the wanted inequality follows from the definition of $J^*(u(.))$: $J^*(u(.)) \leq J(u^*_{\varepsilon,N}(.))$. We need to prove the other side.

For every admissible control u(.), we have

$$J(u(.)) \geq J_{[N]}(u(.)) - \omega_N(t) \geq J^*_{[N]} - \omega_N(T) \geq J_{[N]}(u^*_{\varepsilon,N}(.)) - \varepsilon - \omega_N(t)$$

--- (2.22)

Also, as a consequence of (2.19), we have

$$J_{[N]}(u^*_{\varepsilon,N}(.)) \geq J(u^*_{\varepsilon,N}(.)) - \omega_N(T)$$

--- (2.23)

It follows from (2.22) and (2.23) that, for every admissible control function u(.), we have

$$J(u(.)) \geq J(u^*_{\varepsilon,N}(.)) - \varepsilon - 2\omega_N(T)$$

--- (2.24)

and, consequently,

$$J^*(u(.)) \geq J(u^*_{\varepsilon,N}(.)) - \varepsilon - 2\omega_N(T)$$

--- (2.25)

///

## 3. Parametrization and dynamic programming.

The method of dynamic programing for ODE systems produces a value function (Hamilton- Jacobi function), which representes the infimum of a parametrized cost, where the parametrization is done by initial time and initial state of the controlled system. An $\varepsilon$-suboptimal control can generally be synthesized from the Hamilton-Jacobi function. At the same time, the Hamilton-Jacobi function is important in its own right, as it gives information about the best possible performance of the system under a variety of initial conditions. For this reason, we have to examine the performance, under the proposed approximation acheme, of the family of parametrized cost functionals for the Volterra control problem.

We consider the parametrized problems and the corresponding parametrized cost functionals:

$x_{t,\xi}(\tau)$, $\tau \geq t$ is the solution of the original state dynamics

$$x_{t,\xi}(\tau) = x_0(\tau) + \int_0^\tau f(\tau, s, x_{t,\xi}(s), u(s))\,ds$$

--- (3.1)

subject to the conditions

$$\int_0^t \frac{\partial^{n-1}}{\partial t^{n-1}} f(t, s, x_{t,\xi}(s), u(s))\,ds = \xi_n \quad \forall n \in \mathbb{N}$$

--- (3.2)

$x_{[N],t,\xi_{[N]}}(\tau)$, $\tau \geq t$ is the solution of

$$x_{[N],t,\xi_{[N]}}(\tau) = x_0(\tau) + \int_0^\tau f_{[N]}(\tau, s, x_{[N],t,\xi_{[N]}}(s), u(s))\,ds$$

--- (3.3)

subject to the condition

$$\int_0^t \frac{\partial^{n-1}}{\partial t^{n-1}} f(t, s, x_{[N],t,\xi_{[N]}}(s), u(s))\,ds = \xi_n \quad \forall n = 1, 2, ..., N+1$$

Here, the function $f_{[N]}$ has the same meaning as in section 2. The quantity $\xi_{[N]}$ of (3.3) is the truncation of the $\xi$ of (3.1) to its first N+1 components; in other words, eqns. (3.2) and (3.4) are connected by using the same control u(.) and the same $\xi$.

We shall find another equation that is satisfied by $x_{t,\xi}$. We obtain, from (3.1),

$$x_{t,\xi}(\tau) = x_0(\tau) + \int_0^t f(\tau,s,x_{t,\xi}(s),u(s))\,ds + \int_t^\tau f(\tau,s,x_{t,\xi}(s),u(s))\,ds =$$

$$= x_0(\tau) + \int_0^t \sum_{n=1}^\infty \frac{1}{(n-1)!}\left(\frac{\partial^{n-1}}{\partial t^{n-1}} f(t,s,x_{t,\xi}(s),u(s))\right)(\tau-t)^{n-1}\,ds +$$

$$+ \int_t^\tau f(\tau,s,x_{t,\xi}(s),u(s))\,ds =$$

$$= x_0(\tau) + \sum_{n=1}^\infty \frac{1}{(n-1)!}\xi_n(\tau-t)^{n-1} + \int_t^\tau f(\tau,s,x_{t,\xi}(s),u(s))\,ds$$

--- (3.5)

We set

$$X_0(\tau,t,\xi) := x_0(\tau) + \sum_{n=1}^\infty \frac{1}{(n-1)!}\xi_n(\tau-t)^{n-1}$$

--- (3.6)

Then (3.5) can be rewritten as

$$x_{t,\xi}(\tau) = X_0(\tau,t,\xi) + \int_t^\tau f(\tau,s,x_{t,\xi}(s),u(s))\,ds$$

--- (3.7)

In an analogous manner, we define

$$X_{[N],0}(\tau, t, \xi_{[N]}) := x_0(\tau) + \sum_{n=1}^{N+1} \frac{1}{(n-1)!} \xi_n (\tau - t)^{n-1}$$

--- (3.8)

Then, by the same technique as above, we find that $x_{[N],t,\xi_{[N]}}$ solves the integral equation

$$x_{[N],t,\xi_{[N]}}(\tau) = X_{[N],0}(\tau, t, \xi) + \int_t^\tau f_{[N]}(\tau, s, x_{t,\xi_{[N]}}(s), u(s)) ds$$

--- (3.9)

We shall prove:

<u>Theorem 3.1.</u>  Under the conditions of section 2, the difference $|x_{t,\xi}(\tau) - x_{[N],t,\xi}(\tau)|$ goes to 0, as $N \to \infty$, uniformly for $\tau \in [t, T]$, and the rate of convergence is estimated by the function $\zeta_N(\tau, t, \xi)$ defined in (3.13) below.

<u>Proof:</u>  We have

$$|x_{t,\xi}(\tau) - x_{[N],t,\xi_{[N]}}(\tau)| \leq |X_0(\tau, t, \xi) - X_{[N],0}(\tau, t, \xi_{[N]})| +$$

$$+ \int_t^\tau |f(\tau, s, x_{t,\xi}(s), u(s)) - f_{[N]}(\tau, s, x_{[N],t,\xi_{[N]}}(s), u(s))| ds \leq$$

$$\leq |X_0(\tau, t, \xi) - X_{[N],0}(\tau, t, \xi_{[N]})| +$$

$$+ \int_t^\tau |f(\tau, s, x_{t,\xi}(s), u(s)) - f(\tau, s, x_{[N],t,\xi_{[N]}}(s), u(s))| ds +$$

$$+ \int_t^\tau |f(\tau, s, x_{[N],t,\xi_{[N]}}(s), u(s)) - f_{[N]}(\tau, s, x_{[N],t,\xi_{[N]}}(s), u(s))| ds \leq$$

$$\leq \sum_{n=N+2}^\infty \frac{|\xi_n|}{(n-1)!} (\tau - t)^{n-1} + L_0 \int_t^\tau |x_{t,\xi}(s) - x_{[N],t,\xi_{[N]}}(s)| ds +$$

$$+ \sum_{n=N+2}^\infty \frac{(\tau - t)^{n-1}}{(n-1)!} \int_t^\tau \left| \frac{\partial^{n-1}}{\partial t^{n-1}} f(t, s, x_{[N],t,\xi_{[N]}}(s), u(s)) \right| ds$$

--- (3.10)

We set

$$\varphi_N(\tau,t,\xi) := \sum_{n=N+2}^{\infty} \frac{|\xi_n|}{(n-1)!}(\tau-t)^{n-1} +$$

$$+ \sum_{n=N+2}^{\infty} \frac{(\tau-t)^{n-1}}{(n-1)!} \int_t^\tau \left|\frac{\partial^{n-1}}{\partial t^{n-1}} f(t,s,x_{[N],t,\xi_{[N]}}(s),u(s))\right| ds$$

--- (3.11)

Then, under the conditions of section 2, we have

$$\lim_{N\to\infty} \varphi_N(\tau,t,\xi) = 0 \quad \text{uniformly for } \tau \in [t,T]$$

--- (3.12)

Let $\zeta_N(\tau,t,\xi)$ denote the solution of the Volterra integral equation

$$\zeta_N(\tau,t,\xi) = \varphi_N(\tau,t,\xi) + L_0 \int_t^\tau \zeta_N(s,t,\xi) ds$$

--- (3.13)

Then:

(1) $\zeta_N(\tau,t,\xi) = \varphi_N(\tau,t,\xi) + L_0 \int_t^\tau \exp(L_0(\tau-s))\varphi_N(s,t,\xi) ds$.

(2) $\lim_{N\to\infty} \zeta_N(\tau,t,\xi) = 0$ uniformly for $\tau \in [t,T]$ and $u(.) \in \mathbf{U}_{ad}$.

(3) $|x_{t,\xi}(\tau) - x_{[N],t,\xi_{[N]}}(\tau)| \leq \zeta_N(\tau,t,\xi) \quad \forall \tau \in [t,T]$.

The 3 conclusion above are tantamount to the assertion of the theorem. ///

Next, we prove:

Theorem 3.2. Under the conditions of section 2, the sequence $\{V_{[N]}(t,\xi): N = 1,2,...\}$ of solutions of the finite-dimensional dynamic programming equations, assuming existence of a suitable of a solution in a suitable sense (e.g. the so-called "viscosity solutions") so that the solution indeed represents the optimal parametrized cost, converges to the optimal cost for the original Volterra control problem, and the rate of convergence can be estimated by the function $\psi_N(T,t,\xi)$ defined by

$$\psi_N(T,t,\xi) := L_F \int_t^T \zeta_N(\tau,t,\xi)\,d\tau + L_{F_0}\zeta_N(T,t,\xi)$$

--- (3.14)

where $L_F, L_{F_0}$ are Lipschitz constants of $F$, $F_0$, respectively.

Proof: We have

$$|J_{t,\xi}(u(.)) - J_{[N],t,\xi_{[N]}}(u(.))| \leq \psi_N(T,t,\xi)$$

--- (3.15)

from which it follows that

$$|J^*_{t,\xi} - V_{[N]}(t,\xi)| \leq \psi_N(T,t,\xi)$$

--- (3.16)

according to the following argument:

we have $J^*_{t,\xi} \leq J_{t,\xi}(u(.)) \;\; \forall u(.) \in U_{ad}$, and
$\forall \varepsilon > 0 \; \exists u_\varepsilon(.) \in U_{ad}: V_{[N]}(t,\xi) \geq J_{[N],t,\xi_{[N]}}(u_\varepsilon(.)) - \varepsilon$, thus
$J^*_{t,\xi} - V_{[N]}(t,\xi_{[N]}) \leq J_{t,\xi}(u_\varepsilon(.)) - J_{[N],t,\xi_{[N]}}(u_\varepsilon(.)) + \varepsilon \leq \psi_N(T,t,\xi) + \varepsilon$, and, by taking the limit as $\varepsilon \to 0^+$, we conclude that $J^*_{t,\xi} - V_{[N]}(t,\xi_{[N]}) \leq \psi_N(T,t,\xi)$. In a similar way, we obtain $V_{[N]}(t,\xi_{[N]}) - J^*_{t,\xi} \leq \psi_N(T,t,\xi)$. Thus $|J^*_{t,\xi} - V_{[N]}(t,\xi)| \leq \psi_N(T,t,\xi)$. ///

## 4. Connection between the infinite-dimensional ODE problem and the adjoint integral equations for Volterra control.

The optimal control problem for the infinite-dimensional ODE system (2.3) has, formally, an adjoint infinite-dimensional system of Hamiltonian equations.

The <u>question</u> we address in this section is: what is the connection, if any, between the infinite-dimensional Hamiltonian equations for the system (2.5) and the adjoint integral equations that form the counterpart of Hamiltonian equations for controlled Volterra integral equations?

This question is important for the following reasons:

(1) The co-states of an ODE control problem make up the space-gradient of the Hamilton-Jacobi function. If an approximation to the Hamilton-Jacobi function has been found by using the approximate dynamic programming equations, one way to check its accuracy would be to see if the space-gradient of that solution satisfies (approximately) the adjoint equations. For this purpose, it is useful to use the adjoint equations for the original Volterra control problem, because of the lower dimensionality of the co-state for th original problem compared to the dimensionality of the co-state for the approximating ODE system. Consequently, we need to know the connection between the space-gradient of the Hamilton-Jacobi function and the adjoint integral equations for th original Volterra control problem.

(2) It may be possible to devise iterative methods using simultaneously both the dynamic programming equations for the ODE system (2. ) and the adjoint integral equations for the original Volterra control problem. Again, in that case, we need to know the connection between the two types of co-states, the co-state for the ODE system and the co-state for the original Volterra control problem.

(3) It is mathematically significant to know what is the connection between two different approaches to solving the same optimal control problem.

The <u>answer</u> to this question is: under the condition of continuously differentiable solutions of the dynamic programming equation (2.5), the co-states of the infinite-dimensional ODE control problem are given by

$$\lambda_j(t) = \int_t^T \left\{ \frac{(s-t)^j}{j!} \psi(s) + \frac{(s-t)^{j-1}}{(j-1)!} \frac{\partial}{\partial x_1} F_0(x_0(s) + x_1(s)) \right\} ds.$$

As in the case of similar questions in the theory of optimal control of ODE systems, the connection between dynamic programming and Hamiltonian co-states is established under conditions of existence of classical solutions of the dynamic programming equations. Other possible refinements, using concepts of generalized derivatives, are beyond our present scope.

We write the cost functional J in the form

$$J = \int_0^T F(t, x_0(t) + x_1(t), u(t)) dt + F_0(x_0(T) + x_1(T))$$

--- (4.1)

Since $x_0(.)$ is a fixed function, we may take $x_1(.)$ as the unknown function (and thus $x_1$ as the state variable for the adjoint integral equations for the original Volterra control problem. By standard results on optimal control of Volterra integral equations [M, S, V], the adjoint Volterra integral equation for the co-state $\psi$ is

$$\psi(t) = \frac{\partial}{\partial x_1} F(t, x_0(t) + x_1(t), u(t)) +$$
$$+ \left( \frac{\partial}{\partial x_1} F_0(x_0(T) + x_1(T)) \right) \left( \frac{\partial}{\partial x_1} f(T, t, x_0(t) + x_1(t), u(t)) \right) +$$
$$+ \int_t^T \psi(s) \frac{\partial}{\partial x_1} f(s, t, x_0(t) + x_1(t), u(t)) ds$$

--- (4.2)

We denote by $(\lambda_1(t), \lambda_2(t), ...)$ the co-state for the problem of minimizing J, as given by (4.1), under the state dynamics (2.3). Then the infinite-dimensional Hamiltonian equations are

$$\frac{d\lambda_1(t)}{dt} + \sum_{i=1}^\infty \lambda_i(t) \left( \frac{\partial^{i-1}}{\partial t^{i-1}} f(t, s, x_0(s) + x_1(s), u(s))|_{s=t} \right) + \frac{\partial}{\partial x_1} F(t, x_0(t) + x_1(t), u(t)) = 0;$$

$$\frac{d\lambda_j(t)}{dt} + \lambda_{j-1}(t), \forall j \geq 2; \lambda_1(T) = \frac{\partial}{\partial x_1} F_0(x_0(T) + x_1(T)); \lambda_j(T) = 0, \forall j \geq 2$$

--- (4.3)

The existence of such co-states follows from the assumption of existence of classical solutions of the dynamic programming equations, but it can also be justified by a direct proof of validity of Hamiltonian equations, in a way similar (ut not identical) to the standard case of ODE controlled systems. For our present purposes, a direct derivation of Hamiltonian equations is not necessary.

Now, we define a collection of functions $(\mu_1(t), \mu_2(t), \ldots)$ in the following way:

$$\mu_1(t) := \int_t^T \psi(s)\,ds + \frac{\partial}{\partial x_1} F_0(x_0(t) + x_1(t));\ \mu_j(t) := \int_t^T \mu_{j-1}(s)\,ds,\ \forall j \geq 2$$

--- (4.4)

Then we have

$$\int_t^T \psi(s) \frac{\partial}{\partial x_1} f(s, t, x_0(t) + x_1(t), u(t))\,ds = -\int_t^T \frac{d\mu_1(s)}{ds} \frac{\partial}{\partial x_1} f(s, t, x_0(t) + x_1(t), u(t))\,ds =$$

$$= -\mu_1(t) \left( \frac{\partial}{\partial x_1} f(s, t, x_0(t) + x_1(t), u(t)) \right)\bigg|_{s=t} -$$

$$- \left( \frac{\partial}{\partial x_1} F_0(x_0(T) + x_1(T)) \right) \left( \frac{\partial}{\partial x_1} f(T, t, x_0(t) + x_1(t), u(t)) \right) +$$

$$+ \int_t^T \mu_1(s) \frac{\partial^2}{\partial s \partial x_1} f(s, t, x_0(t) + x_1(t), u(t))\,ds$$

--- (4.5)

and, inductively, $\forall j \geq 2$,

$$\int_t^T \mu_{j-1}(s) \frac{\partial^j}{\partial^{j-1} s \partial x_1} f(s, t, x_0(t) + x_1(t), u(t))\,ds =$$

$$= -\int_t^T \frac{d\mu_j(s)}{ds} \frac{\partial^j}{\partial^{j-1} s \partial x_1} f(s, t, x_0(t) + x_1(t), u(t))\,ds =$$

$$= -\mu_j(t) \left( \frac{\partial^j}{\partial s^{j-1} \partial x_1} f(s, t, x_0(t) + x_1(t), u(t)) \right)\bigg|_{s=t} +$$

$$+ \int_t^T \mu_j(s) \frac{\partial^{j+1}}{\partial^j s \partial x_1} f(s, t, x_0(t) + x_1(t), u(t))\,ds$$

--- (4.6)

Under the condition

$$\lim_{j\to\infty} \int_t^T \mu_j(s) \frac{\partial^{j+1}}{\partial^j s \partial x_1} f(s,t,x_0(t)+x_1(t),u(t))\,ds = 0 \text{ uniformly for } t \in [0,T] \text{ and } u(.) \in \mathbf{U}_{ad}$$

--- (4.7)

we find, from (4.5) and (4.6), that

$$\int_t^T \psi(s) \frac{\partial}{\partial x_1} f(s,t,x_0(t)+x_1(t),u(t))\,ds = -\mu_1(t)\left(\frac{\partial}{\partial x_1} f(s,t,x_0(t)+x_1(t),u(t))\right)\Big|_{s=t}$$
$$- \sum_{j=2}^{\infty} \mu_j(t)\left(\frac{\partial^j}{\partial s^{j-1}\partial x_1} f(s,t,x_0(t)+x_1(t),u(t))\right)\Big|_{s=t}$$
$$- \left(\frac{\partial}{\partial x_1} F_0(x_0(T)+x_1(T))\right)\left(\frac{\partial}{\partial x_1} f(T,t,x_0(t)+x_1(t),u(t))\right)$$

--- (4.8)

We also note the following consequences of (4.4):

$$\psi(t) = -\frac{d\mu_1(t)}{dt};\; \mu_{j-1}(t) = -\frac{d\mu_j(t)}{dt}\;\; \forall j \geq 2;$$

$$\mu_1(T) = \frac{\partial}{\partial x_1} F_0(x_0(T)+x_1(T));\; \mu_j(T) = 0 \;\forall j \geq 2$$

--- (4.9)

In view of (4.8) and (4.9), we can transform (4.2) into

$$\frac{d\mu_1(t)}{dt} + \mu_1(t)\left(\frac{\partial}{\partial x_1}f(s,t,x_0(t)+x_1(t),u(t))\right)\bigg|_{s=t} +$$

$$+ \sum_{j=2}^{\infty} \mu_j(t)\left(\frac{\partial^j}{\partial s^{j-1}\partial x_1}f(s,t,x_0(t)+x_1(t),u(t))\right)\bigg|_{s=t} + \frac{\partial}{\partial x_1}F(t,x_0(t)+x_1(t),u(t)) = 0;$$

$$\mu_{j-1}(t) = -\frac{d\mu_j(t)}{dt} \ \forall j \geq 2; \ \mu_1(T) = \frac{\partial}{\partial x_1}F_0(x_0(T)+x_1(T)); \ \mu_j(T) = 0 \ \forall j \geq 2$$

--- (4.10)

Since

$$\left(\frac{\partial^j}{\partial s^{j-1}\partial x_1}f(s,t,x_0(t)+x_1(t),u(t))\right)\bigg|_{s=t} = \left(\frac{\partial^j}{\partial t^{j-1}\partial x_1}f(t,s,x_0(s)+x_1(s),u(s))\right)\bigg|_{s=t}$$

--- (4.11)

we see that (4.10) is identical to (4.3) with $\mu_j$, $j=1,2,...$ in lieu of $\lambda_j$, $j=1,2,...$ .

Under the condition that (4.3) possesses a unique solution, we conclude that

$$\lambda_1(t) = \int_t^T \psi(s)ds + \frac{\partial}{\partial x_1}F_0(x_0(t)+x_1(t)); \ \lambda_j(t) = \int_t^T \lambda_{j-1}(s)ds, \ \forall j \geq 2$$

--- (4.12)

Thus

$$\lambda_j(t) = \int_t^T \int_{s_1}^T \cdots \int_{s_{j-1}}^T \psi(s_j)ds_j \cdots ds_2\, ds_1 +$$

$$+ \int_t^T \int_{s_1}^T \cdots \int_{s_{j-2}}^T \left(\frac{\partial}{\partial x_1}F_0(x_0(s_{j-1})+x_1(s_{j-1}))\right)ds_{j-1}\cdots ds_2\, dt$$

--- (4.13)

By changing the order of integration in the formula for $\lambda_j(t)$ and carrying out the integrations with respect to all other variables except the variable appearing inside $\psi$, we find

$$\lambda_j(t) = \int_t^T \left\{ \frac{(s-t)^j}{j!} \psi(s) + \frac{(s-t)^{j-1}}{(j-1)!} \frac{\partial}{\partial x_1} F_0(x_0(s) + x_1(s)) \right\} ds$$

--- (4.14)

In (4.14), for the case j=1, we define $\lambda_1(T)$ by continuity as $t \to T^-$, i.e.

$$\lambda_1(T) = \frac{\partial}{\partial x_1} F_0(x_0(T) + x_1(T)).$$

Also, it follows either from (4.12) or from (4.13) that

$$\psi_j(t) = (-1)^j \frac{d^j \lambda_j(t)}{dt^j}$$

--- (4.15)

Eqns. (4.14) and (4.15) provide the connection between the co-state for the infinite-dimensional ODE control problem and the co-state for the original Volterra control problem.